\documentclass{article}
\usepackage{amssymb,amsthm,amsfonts,amsmath}

\usepackage{hyperref}

\newcommand{\eq}{\begin{equation}}
\newcommand{\en}{\end{equation}}

\newtheorem{theorem}{\large Theorem}%[section]

\newtheorem{lemma}[theorem]{\large Lemma}

\begin{document}
\title{Moment problems and boundaries of number triangles
}

\author{Alexander Gnedin\thanks{Utrecht University; e-mail gnedin@math.uu.nl}
\hspace{.2cm}
and 
\hspace{.2cm}
Jim Pitman\thanks{University of California, Berkeley; e-mail pitman@stat.Berkeley.EDU} 
}
\date{
\today
\\
}
\maketitle
\begin{abstract}
The boundary problem for graphs like Pascal's but with general multiplicities of edges 
is related to a `backward' problem of moments of the Hausdorff type.\\

 
\end{abstract}

\section{The extreme boundary}

Let $T_n:=\{(n,0),(n,1),\ldots,(n,n)\}$ and $T:=\cup_{n=0}^\infty T_{n}$.
We endow $T$  with the structure of a directed graph in which every node $(n,k)$ has two outgoing edges $(n,k)\to (n+1,k)$ and 
$(n,k)\to (n+1,k+1)$ with {\em multiplicities}  $\ell_{nk}$ and $r_{nk}$ (respectively), where
$\{\ell_{nk};~  (n,k)\in T\}$ 
and $\{r_{nk};~  (n,k)\in T\}$ are given
triangular arrays with (strictly) positive entries.
A classical example is the Pascal graph with unit multiplicities $\ell_{nk}=r_{nk}=1$.

\par Let $\cal V$ be the set of nonnegative solutions 
$V=\{V_{nk};~  (n,k)\in T\}$
to the backward recursion
\begin{equation}\label{Vrec}
V_{nk}=\ell_{nk} V_{n+1,k}+ r_{nk} V_{n+1,k},~~~~~~~(n,k)\in T
\end{equation}
with normalisation $V_{00}=1$. The set $\cal V$  is convex and compact
in the product topology of functions on $T$.
By some general theory in Dynkin (1978)
the {\em extreme boundary} ${\rm ext} {\cal V}$,
comprised of  indecomposable elements of $\cal V$,
is a Borel set. 
Moreover,
$\cal V$ is a Choquet simplex,
 meaning that 
each $V\in {\cal V}$ has a unique representation as convex combination 
\begin{equation}\label{VU}
V=\int_{{\rm ext}{\cal V}} U\,\mu({\rm d}U)
\end{equation}
with some probability measure $\mu$ supported by 
${\rm ext} {\cal V}$. The {\em boundary problem} for the graph $T$ is
to find some explicit description of the set of extremes, meaning, if possible, 
a  simple parametrisation of   
${\rm ext} {\cal V}$ along with the kernel that is implicit in (\ref{VU}).

 \par The recursion  (\ref{Vrec}) for the Pascal graph appeared in the work of
Hausdorff on summation methods (1921, p. 78) and the `small' problem of moments
on $[0,1]$.
In this case the bivariate array $V$ is completely determined by $V_{\bullet,0}$ according to the rule 
$V_{\bullet+k,k}=\nabla^k (V_{\bullet,0})$,
where $\bullet$ stands for the variable $n$, and $\nabla^k$ is the $k$th iterate of the difference operator $\nabla (U_\bullet):=U_\bullet-U_{\bullet+1}$.
The condition $V\geq 0$ means that $V_{\bullet, 0}$ is completely 
monotone, hence by Hausdorff's theorem $V_{\bullet, 0}$ is a sequence of moments 
$$V_{n0}=\int_{[0,1]} x^n \mu({\rm d}x)$$
of some probability measure $\mu$.
That is to say, the set  of extremes ${\rm ext}{\cal V}$ can be identified with the unit interval, 
and the extremes  have the form $V_{nk}(x)=x^{n-k}(1-x)^k$
for $x\in [0,1]$; in particular, $V_{n0}=x^n$.

\par For general multiplicities  the recursion (\ref{Vrec}) is equivalent to
$$V_{\bullet, k}=\nabla_k(\cdots(\nabla_1(V_{\bullet,0}))\cdots),$$
where $\nabla_k(U_{\bullet})=(U_{\bullet}-\ell_{\bullet,k}U_{\bullet+1})/r_{\bullet, k}$
is a generalised difference operator. 
%In view of this,
%the condition $V\geq 0$ is  analogous to the complete monotonicity of sequence  $V_{\bullet,0}$, 
By analogy with Hausdorff's criterion, 
the question about positivity of the generalised iterated differences of $V_{\bullet,0}$
may be regarded as a `backward'  problem of moments. A direct problem of moments of the Hausdorff 
type appears when we determine the $V_{n0}$'s for extreme solutions as 
 functions on the boundary, and consider 
the integral representation of the generic $V_{n0}$ in the form (\ref{VU}). 

%the integral representation  of $V$ in the form  
%(\ref{VU}) may be regarded as 
%a generalised problem of moments, understood as a criterion of positivity of the generalised differences of 
%$V_{\bullet,0}$.  

\par A bivariate array $V\in {\cal V}$ 
could be also computed by suitable differencing 
the diagonal sequence $(V_{nn})$, but 
this leads to the same type of the moment problem by virtue of the transposition of $T$  which exchanges  
the multiplicities $\ell_{nk}$'s with $r_{n,n-k}$'s.

\par A special feature of $T$, as compared with more complicated graphs like  Young's lattice 
(see Kerov (2003), Borodin and Olshanski (2000)),
is a natural total order on the extreme boundary.
In this note we extend the argument of 
 Gnedin and Pitman (2006) to show that the total order allows ${\rm ext}{\cal V}$ to be 
embedded into $[0,1]$.  We shall also survey the 
connection of the boundary problem with asymptotic properties
of some classical arrays of combinatorial numbers.

\section{Markov chain approach}

The {\em weight} of a path in $T$ joining the root $(0,0)$ and some other node $(n,k)$ is defined as the product of multiplicities of edges along the path
(for instance the weight of $(0,0)\to(1,0)\to(2,1)$ is $\ell_{00}r_{10}$).
The {\em dimension} $D_{nk}$ of $(n,k)\in T$ is defined to be the sum of weights of all paths from $(0,0)$ to $(n,k)$. 
The dimensions are computable from the  forward recursion
\begin{equation}\label{Drec}
D_{nk}=r_{n-1,k-1} D_{n-1,k-1}+ \ell_{n-1,k} D_{n-1,k},~~~~~~~(n,k)\in T,
\end{equation}
(where the first term in the right-hand side is absent for $k=0$ and the second term is absent for $k=n$),
 with the initial condition $D_{00}=1$.
The {\em number triangle} associated with $T$ is the array $\{D_{nk};~  (n,k)\in T\}$.

\par 
Each $V\in {\cal V}$ determines the law ${\mathbb P}_V$ of a inhomogeneous Markov chain $K_\bullet$ whose 
 backward transition probabilities    for $0\leq k\leq n,~n>0$ are   
\begin{equation}\label{backtr}
{\mathbb P}_V(K_{n-1}=j\,|\, K_n=k)= {D_{n-1,j}\over D_{nk}}(\ell_{n-1,j}\delta_{jk}+r_{n-1,j}\delta_{j,k-1}), 
\end{equation}
and whose distribution at time $n$ is
${\mathbb P}_V(K_n=k)=D_{nk} V_{nk}.$
It is important that the probabilities (\ref{backtr}) are determined solely by the multiplicities of edges and do not depend on $V$.
Hence  $\cal V$ is in essence a class of distributions for Markov chains on $T$ with given backward transition probabilities.

\par  For each fixed  integer $\nu$ and $0\leq \varkappa\leq \nu$  
let $V^{\nu\varkappa}$ be the function on $T$ which satisfies the recursion (\ref{Vrec}) for 
$n< \nu$, satisfies
$V^{\nu\varkappa}_{\nu k}=\delta_{k\varkappa}$, 
and equals $0$ on $\cup_{n>\nu}T_n$. 
Such $V^{\nu \varkappa}$ determines the probability law of a finite Markov chain $(K_0,\ldots,K_\nu)$
conditioned on $K_\nu=\varkappa$.

\par We define the {\em sequential boundary} $\partial_\uparrow{\cal V}$ to be the set of elements of ${\cal V}$ representable as limits 
$V=\lim_{\nu\to\infty} V^{\nu,\varkappa(\nu)}$
taken along  infinite paths   $\{\varkappa(\nu);~\nu=0,1,\ldots\}$  in $T$.
The sequential boundary $\partial_\uparrow{\cal V}$ may be smaller than
the  set of all accumulation points for $\{V^{\nu\varkappa}; (\nu,\varkappa)\in T\}$
(the {\em Martin boundary}), but it is large enough to cover ${\rm ext}{\cal V}$, as is seen from 
the following  lemma, which is a variation on the theme of sufficiency (see Diaconis and Freedman (1984)).

\begin{lemma} If $V\in {\rm ext}{\cal V}$ then the random functions
\begin{equation}\label{defVX}
V^{\nu,K_\nu}:=\sum_{\varkappa=0}^{\nu} 1_{\{K_\nu= \varkappa\}}   V^{\nu\varkappa}
\end{equation}
satisfy
$V^{\nu,K_\nu}\to V$, as $\nu\to\infty$, ${\mathbb P}_V$-almost surely.
\end{lemma}
\begin{proof}
 Let ${\cal F}_\nu$ be the sigma-algebra generated by $K_{\nu}, K_{\nu+1}, \ldots$, and ${\cal F}_\infty=\cap {\cal F}_\nu$.
Let ${\mathbb P}_V$ correspond to some extreme $V$.
Choose any $(n,k)$ and consider random variables 
$$V^{\nu,K_\nu}_{nk}={\mathbb P}_V(K_n=(n,k)\,|\, K_\nu)/D_{nk}={\mathbb P}_V(K_n=(n,k)\,|\, {\cal F}_\nu)/D_{nk},~~~\nu\geq n$$
where the first equality follows from the definition (\ref{defVX}), and the second equality is a consequence of the Markov 
property. 
Applying  Doob's reversed martingale convergence theorem to the conditional expectations given ${\cal F}_1\supset {\cal F}_2\supset \ldots$
we obtain
$$V^{\nu, K_\nu}_{nk}\to {\mathbb P}_V(K_n=(n,k)\,|\, {\cal F}_\infty)/D_{nk}~~~~~{\mathbb P}_V{\rm-a.s.}$$
The assumption $V\in{\rm ext} {\cal V}$ implies that $\cal F_\infty$ is trivial, hence
$${\mathbb P}_V(K_n=(n,k)\,|\, {\cal F}_\infty)={\mathbb P}_V(K_n=(n,k))=V_{nk}D_{nk}.$$
%Since $(n,k)$ is arbitrary, we are done.
\end{proof}
\vskip0.3cm
\par Thus ${\rm ext}{\cal V}\subset\partial_\uparrow{\cal V}$ (in general the inclusion is strict). To state this conclusion in analytical terms,
define the weight of a path in $T$ connecting two nodes $(n,k)$ and $(\nu,\varkappa)$ as the product of multiplicities along the path,
and define the {\em extended dimension} $D^{\nu\varkappa}_{nk}$ as the sum of weights over all such paths
(so   that $D^{\nu\varkappa}_{00}=D_{\nu\varkappa}$). We then have a  fundamental relation
\begin{equation}\label{reldim}
V_{nk}^{\nu\varkappa}={D_{nk}^{\nu\varkappa}\over D^{\nu\varkappa}},
\end{equation}
which connects the boundary problem with asymptotic properties of $T$.
Specifically, the convergence of $V^{\nu,\varkappa(\nu)}$ amounts to the convergence of these ratios for all $(n,k)\in T$ along the path
(in fact, it is enough to focus on $V_{\bullet,0}$).

\section{Order}

A special feature of $T$ which yields the order  is that the only possible increments 
of the variable $k$ along any path  are $0$ and $1$. 
The next lemma appeared in Gnedin and Pitman (2006) with a different proof.

\begin{lemma}\label{L2} For $\nu> n$ fixed, $V_{n0}^{\nu\varkappa}$ is nonincreasing in $\varkappa$.
\end{lemma}
\begin{proof} Choose  $0\leq\varkappa<\varkappa'\leq \nu$ and consider two Markov chains $K_\bullet, K_\bullet'$ 
which run in reverse time $n=\nu,\nu-1,\ldots,0$ according to (\ref{backtr}) and start with $K_\nu=\varkappa, K_\nu'=\varkappa'$. 
Suppose the chain $K_\bullet'$ jumps independently of $K_\bullet$ 
as long as they are 
in distinct states, and suppose that  $K_\bullet'$ is coupled with $K_\bullet$  
at some random time $0\leq \tau<\nu$ when the states become the same. In the reverse time, only transitions $k\to k,~ k\to k-1$ for $k>0$
and $0\to 0$ are possible, hence we always have $K_n'\geq K_n$. Therefore the event $K'_n=0$ occurs exactly when $K_n=0$ and $\tau\geq n$, which implies 
$${\mathbb P}(K_n=0\,|\, K_\nu=\varkappa)\geq {\mathbb P}(K_n=0\,|\, K_\nu=\varkappa').$$
\end{proof}
\vskip0.3cm

\par A minor modification of the above argument shows that if $K_\nu$ under ${\mathbb P}_V$ is strictly 
stochastically smaller than $K_\nu$ under some other ${\mathbb P}_{V'}$,
then the same relation holds true for every $n\leq \nu$.

\par We focus now on $V_{10}$.  Suppose $V\in \partial_\uparrow{\cal V}$ is induced, via (\ref{reldim}), by 
some infinite path $\{\varkappa(\nu);~\nu=0,1,\ldots\}$, and  $V'\in \partial_\uparrow{\cal V}$ is induced by 
some other path $\{\varkappa'(\nu);~\nu=0,1,\ldots\}$. If $\varkappa(\nu)=\varkappa'(\nu)$ for infinitely many $\nu$ then, of course,
$V=V'$. If $\varkappa(\nu)<\varkappa'(\nu)$ for infinitely many $\nu$ and 
$\varkappa(\nu)>\varkappa'(\nu)$ for infinitely many $\nu$ then by Lemma \ref{L2} we have $V_{\bullet, 0}=V_{\bullet, 0}'$ and $V=V'$. 
Thus $V\neq V'$ can only occur if the same strict inequality holds for all sufficiently large $\nu$.
To be definite, let $\varkappa(\nu)< \varkappa'(\nu)$ for all large enough $\nu$,
but then $V\neq V'$ implies that 
$K_n$ under ${\mathbb P}_V$ is strictly 
stochastically smaller than $K_n$ under ${\mathbb P}_{V'}$
for all $n>0$, in particular this holds for $n=1$ which means that $V_{10}>V_{10}'$. 
We see that
for $V,V'\in \partial_\uparrow {\cal V}$, the inequality $V_{10}>V_{10}'$ holds if and only if $K_n$ under ${\mathbb P}_V$ is strictly stochastically
smaller than $K_n$ under ${\mathbb P}_{V'}$ for all $n>0$.
This defines a strict order  $\lhd$  on $\partial_\uparrow {\cal V}$.

\begin{lemma} The sequential boundary $\partial_\uparrow {\cal V}$ is compact.
\end{lemma}
\begin{proof} Suppose $V^j\in \partial_\uparrow {\cal V}$ $(j=1,2,\ldots)$ is a  sequence converging to some $V\in{\cal V}$. 
We know that $\cal V$ is a metrisable compactum with some distance function $\rm dist$.
Passing to a subsequence 
we can restrict consideration to the case of
increasing or decreasing sequence, 
so to be definite assume that $V^{j+1}\lhd V^j$ for $j=1,2,\ldots$
Choosing some path $\{\varkappa^j(\nu);~ \nu=0,1,\ldots\}$ which induces $V^j$, the ordering implies that
$\varkappa^j(\nu)\to\infty$ as $\nu\to\infty$
and 
$\varkappa^j(\nu)<\varkappa^{j+1}(\nu)$ for all large enough $\nu$.
As $\nu$ varies, define inductively in $j$ 
a function $\varkappa(\nu)$ which coincides for some $\nu$ with $\varkappa^j(\nu)$.
Specifically,
$\varkappa(\nu)=\varkappa^j(\nu)$ until $\varkappa^{j+1}(\nu)<\varkappa^j(\nu)$ starts to hold 
along with ${\rm dist}(V^{\nu,\varkappa_j},V)<1/j$ and ${\rm dist}(V^{\nu,\varkappa_{j+1}},V)<1/j$, 
then let $\varkappa(\nu)$ decrement
by $1$ until it becomes equal to $\varkappa^{j+1}(\nu)$. 
This defines an infinite path in $T$, for which one can use 
monotonicity to show that
$V^{\nu,\varkappa(\nu)}\to V$.
\end{proof}

\par Recalling that $\ell_{00}V_{10}+r_{00}V_{11}=1$ we obtain:

\begin{theorem}
The function $V\mapsto \ell_{00}V_{10}$ is an ordered  homeomorphism of the sequential boundary $\partial_{\uparrow}{\cal V}$ with order
$\lhd$ into
$[0,1]$ with order $>$.
\end{theorem}

\noindent
Two extreme cases $\ell_{00}V_{01}=0$ and $\ell_{00}V_{01}=1$ correspond to trivial Markov chains
$K_\bullet=0$ and $K_\bullet=\bullet$, respectively.

\section{Discrete or continuous?}

In the situation covered by the following lemma, 
setting $\varkappa(\nu)=m$ (for large $\nu$) for $m=1,2,\ldots$ is  the only way to induce
 nontrivial limits. Then
 ${\rm ext} {\cal V}$ is discrete
 and coincides with the sequential boundary.

\begin{lemma}{\rm (Gnedin and Pitman (2006))}\label{disc}
Suppose for $m=0,1,\ldots$ there are solutions $V(m) \in {\cal V}$ such that 
$V_{nm} (m)\,D_{nm}\to 1$ as $n\to\infty$, 
then each $V(m)$ is extreme and satisfies
$K_n\to m$ ${\mathbb P}_{V(m)}$-{\rm a.s.}.
If also $V_{10}(m)\to 0$ as $m\to\infty$
then $V(m)$ converges to the trivial solution $V(\infty)$ with 
$K_\bullet=\bullet$ ${\mathbb P}_{V(\infty})$-{\rm a.s.},
and in this case 
${\rm ext}{\cal V}=\partial_\uparrow{\cal V}=\{V(0), V(1),\ldots, V(\infty)\}$. 
\end{lemma}

\par In some cases the limits can be obtained 
by setting $\varkappa(\nu)\sim s\, c(\nu)$ with suitable scaling  $c(\nu)\to\infty$ and $s\geq 0$.
Under conditions in the next lemma, ${\rm ext}{\cal V}$ coincides with $\partial_\uparrow{\cal V}$ and is homeomorphic to $[0,1]$.
The scaling determines the order of growth of $K_\bullet$ under ${\mathbb P}_V$'s.

\begin{lemma}\label{cont}{\rm (Gnedin and Pitman (2006))} 
Suppose there is a sequence of positive constants $\{c(\nu); \nu=0,1,\ldots\}$
with $c(\nu)\to\infty$,
and for each $s\in [0,\infty]$ there is
a solution $V(s)\in {\cal V}$ which satisfies 
$K_\nu/c(\nu)\to s$ ${\mathbb P}_{V(s)}$-{\rm a.s.}.
Suppose 
the mapping $s\mapsto V(s)$ is a continuous injection from $[0,\infty]$ to ${\cal V}$ with
$0$ and $\infty$ corresponding to the trivial solutions.
%$${\mathbb P}_{{\cal V}(0)}(K_n=1)=1\,, ~~~{\mathbb P}_{{\cal V}(\infty)}(K_n=n)=1\,.$$
Then a
 path $\{\varkappa(\nu);\,\nu=0,1,\ldots\}$ induces a limit 
if and only if $\varkappa(\nu)/c(\nu)\to s$ for some 
$s\in [0,\infty]$, in which case the limit is $V(s)$.
Moreover, ${\rm ext}{\cal V}=\partial{\cal  V}_\uparrow=\{V(s), s\in [0,\infty]\}$. 
\end{lemma}

\par Minor variations of the above two situations are obtained by 
 transposing  multiplicities $\ell_{nk}\leftrightarrow r_{n,n-k}$.
%(we need to  replace $\varkappa(\nu)$ 
%through $\nu-\varkappa(\nu)$ and other obvious adjustments. 
Still, this does not exhaust all possibilities.
See Kerov (2003) (Section 1.3, Theorem 2) for examples of
boundaries with both discrete and continuous components.

\section{Number triangles}

\noindent
{\bf The Pascal triangle.}
For the Pascal graph the dimensions are $D_{nk}={n\choose k}$ and $D^{\nu\varkappa}_{nk}={\nu-n\choose \varkappa-k}$.
The ratios $V^{\nu,\varkappa(\nu)}_{nk}={\nu-n\choose \varkappa(\nu)-k}/{\nu\choose \varkappa(\nu)}$ converge iff $\varkappa(\nu)/\nu\to x\in [0,1]$,
in which case the limit is $V_{n0}(x)=x^n$.
This identification of extremes is equivalent to de Finetti's theorem (see  Aldous (2003)), since $V\in {\cal V}$ 
determines the law of some infinite sequence of exchangeable Bernoulli trials.
A closely related type of moment problem with a monotonicity constraint
have been discussed recently in Gnedin and Pitman (2007).

\noindent
{\bf The $q$-Pascal triangle.}
This graph  has multiplicities $\ell_{nk}=1, r_{nk}=q^{n-k}$, $(n,k)\in T$,
and
may be seen as a parametric deformation of the Pascal graph. The extreme boundary 
was  found in  Kerov (2003) by an algebraic method and justified by Olshanski (2001)  by the analysis of 
(\ref{reldim}).
 The dimensions are expressible through  $q$-binomial coefficients as
$$D_{nk}={n\choose k}_q,~~~~~D^{\nu,\varkappa}_{nk}=q^{(\varkappa-k)(n-k)}{\nu-n\choose \varkappa-k}_q{\bigg /}{\nu\choose \varkappa}_q.$$
Suppose first that $0<q<1$.
Lemma \ref{disc} is applicable, and all nontrivial extremes are given by 
$$V_{nk}(m)= {q^{(m-k)(n-k)}(1-q)\cdots (1-q^m)\over(1-q)\cdots(1-q^{m-k})}~ 1_{\{0\leq k\leq m\}}, ~~~~m=1,2,\ldots$$  
 In particular, $V_{n0}(m)=q^{mn}$, $(m=0,1,\ldots,\infty$).

\par 
The function $V\mapsto \ell_{00}V_{10}$ identifies the extreme boundary with  $\{q^m,~m=0,1,\ldots,\infty\}$.
The decomposition (\ref{VU}) into extremes corresponds
to a version of Hausdorff's moment problem on $[0,1]$ with kernel $x^n$, 
but subject to 
the constraint that the measure is to be supported by a geometric progression.
That is to say,
a sequence $V_{\bullet,0}$ with $V_{00}=1$ is representable as a mixture
$$V_{n0}=\sum_{m\in \{0,1,\ldots,\infty\}} p_m q^{mn}$$
with some probability distribution $\{p_m;~ m=0,1,\ldots,\infty\}$ if and only if 
$V_{\bullet, k}=\nabla_k(\cdots(\nabla_1(V_{\bullet, 0})\cdots)\geq 0$ for all $k\geq 0$, where
$\nabla_k(U_\bullet)= (U_\bullet -U_{\bullet+1})/q^{\bullet-k}$.
%In particular, $V_{\bullet, 0}$ is completely monotone.

\par In the case $q>1$ the extreme boundary is  $\{1-q^{-m},~m=0,1,\ldots,\infty\}$ (this case 
is reducible to $q<1$ by transposition of $T$ and replacing $q$ by $q^{-1}$). 
The only accumulation point of ${\rm ext}{\cal V}$ for $q<1$ is $0$ and for $q>1$ is $1$.
A phase transition occurs at $q=1$, when the extreme boundary is continuous.

\noindent
{\bf Stirling triangles.}
Let $r_{nk}=1$ and $\ell_{nk}=(n+1)-\alpha(k+1)$ for $-\infty<\alpha<1$. For $\alpha=-\infty$ take
 $\ell_{nk}=k+1$. 
The dimension is
$D_{nk}=\left[\!\! \begin{array}{c} n+1\\k+1 \end{array}\!\!\right]_\alpha$. The notation stands for 
the generalised Stirling numbers  defined as
connection coefficients in 
$$(t)_{n\uparrow}=\sum_{k=1}^n \left[\!\! \begin{array}{c} n\\k \end{array}\!\!\right]_\alpha 
\alpha^n(t/\alpha)_{n\uparrow},$$
(where $\uparrow$ denotes the rising factorial),
with the convention that these are the Stirling numbers of the second kind for $\alpha=-\infty$.
For $\alpha=0$ these are the  signless Stirling numbers of the first kind.

\par For $-\infty\leq\alpha<0$ the extreme boundary is discrete, with
$$V_{n0}(m)={1\over (m|\alpha|+1)_n}~~~{\rm for~~}-\infty<\alpha<0, ~~~~V_{n0}(m)={1\over m^n}~~~{\rm for~~}\alpha=-\infty.$$
These kernels underly a moment problem for measures on the set $\{0,1,\ldots,\infty\}$.

\par A phase transition occurs at $\alpha=0$.
  Lemma \ref{cont} applies with $\varkappa(\nu)\sim s\,\log n$, the extreme  boundary is continuous and the kernel is
$$V_{n0}(s) ={1\over (s+1)_{n\uparrow}}~~~~~s\in[0,\infty].$$ 
This case is closely related to  random permutations, records and Ewens' sampling formula (see Arratia et al (2003)).

\par In the case  $0<\alpha<1$  we should take $\varkappa(\nu)\sim s\, n^\alpha$ to generate the boundary, see 
 Gnedin and Pitman (2006)  for formulas for 
$V_{n0}(s)$ (to adjust the notation in Gnedin and Pitman (2006) to the present setting, one should replace $(n,k)$ by $(n+1,k+1)$).  
This family of solutions is related to Poisson-Kingman partitions, see Gnedin and Pitman (2006)  and references therein.

\par Several  results and (still open) conjectures about boundaries of more general
Stirling graphs, with multiplicities of the form $\ell_{nk}=b_n+a_k, ~r_{nk}=1$, are 
given in  Kerov (2003).

\noindent
{\bf The Eulerian triangle.}
For multiplicities 
$\ell_{nk}=k+1, ~r_{nk}=n-k+1$
the dimension is the Eulerian number $\langle {n+1\atop k}\rangle$
(that counts permutations with a given number of descents).
The boundary problem was solved in Gnedin and Olshanski (2006).
The extreme solutions are given by
$$V_{nk}(m)={1\over (n+1)!} \prod_{i=-k}^{n-k} \left( 1+{i\over m}\right)$$
with $m\in {\mathbb Z}\cup\{\infty\}$. 
Note that the range of $\ell_{00}V_{10}(m)=(m+1)/(2m)$ is symmetric about $1/2$, with $1/2$ being the only accumulation point.
The symmetry of the boundary stems in this case from the invariance of multiplicities under transposition.

\end{document}